\documentclass[11pt]{article}

\def\N{I\!\!N}

\begin{document}

\vspace*{-3.8cm}

\begin{center}
{\Large {\bf {Dreidel Lasts $O(n^2)$ Spins}}} 

\vskip 10 pt

{\large {\em {Thomas Robinson $\mbox{\qquad \qquad}$ Sujith Vijay}}} 

\vskip 10 pt

{\large {\tt {Department of Mathematics}}}, 

\vskip 5 pt

{\large {\tt {Rutgers, the State University of New Jersey}}} 
\end{center}

\vskip 20 pt

\addtolength{\parskip}{-5mm}

{\em {Dreidel}} is a popular game played during the festival of Chanukah.  
Players start with an equal number of tokens, and contribute one token
each to a common pot. They then take turns spinning a four-sided top,
called the dreidel, and depending on the side showing up, the spinner does
one of the following: \\\\ 
\begin{tabular}{lll}
{\em Nisht} (N) & : & Nothing. \\ {\em Ganz} (G) & : & 
Takes all the tokens in the pot. \\ {\em Halb} (H) & : & Takes (the 
smaller) half of the number of tokens in the pot. \\ 
{\em Shtel} (S) & : & Donates one token to the pot. \\
\end{tabular} \\\\

Whenever the pot is empty, all the players {\em ante up}, i.e., donate one
token each to the pot. Players lose, and {\em go home}, when they are
required to donate a token to the pot, but cannot. The last survivor wins.  
The winner also goes home. \\

\addtolength{\parskip}{2mm}

Feinerman [2] and Trachtenberg [3] investigated the fairness of a
simplified model of dreidel. Zeilberger [4] conjectured that the expected
number of spins in a game of dreidel between two players starting with $n$
tokens each is $O(n^2)$. Later, Banderier [1] conjectured that even in a
multi-player game, the expected number of spins before one of the players
goes home is $O(n^2)$.  We show that the expected duration of a game of
dreidel where the players start with $n$ tokens each is $O(n^2)$,
irrespective of the number of players. \\

Let $P_1,P_2,\ldots,P_k$ denote the players, in the order in which they
spin the dreidel. We introduce three variants of the game of dreidel. \\

{\em {Metadreidel}} works like dreidel, except that the players do not
necessarily start with an equal number of tokens. {\em {Slowdel}} also
works like dreidel, except that it is divided into {\em epochs}, and
allows overdraft, so that the players can continue to play with a negative
number of tokens. Define $k$ spins to be a {\em round}. An epoch ends when
the last spin in a round results in a Ganz (for player $P_k$). The ante
up that follows is also part of the same epoch. A player loses if and only
if he or she has a negative number of tokens at the end of an epoch. {\em
{Metaslowdel}} is the slowdel analogue of metadreidel. \\

Clearly, the slowdel analogue of any instance of a game of dreidel or
metadreidel lasts at least as many spins. \\

Consider a metaslowdel game where $P_k$ starts with $W_0$ tokens, $0 \le
W_0 \le k(n-1)$. Let $W_i$ denote the number of tokens $P_k$ has at the
end of the $i^{th}$ epoch. For $i \ge 1$, we define $Y_i = W_i - W_{i-1}$
to be the {\em payoff} of $P_k$ during the $i^{th}$ epoch. Note that
$\{Y_i\}$ is a set of independent and identically distributed random
variables, with $S_m \doteq \sum_{i=1}^m Y_i = W_m - W_0$. Let $\mu \doteq
E(Y_1)$ and $\sigma^2 \doteq Var(Y_1)= E(Y_1^2) - \mu^2$. \\

Let $T=\inf_{j \in \N} \{j: S_j < - W_0 \mbox{ or } S_j > k(n-1) - W_0
\}$, so that $P_k$ goes home at the end of the $T^{th}$ epoch. Observe
that $T$ is a stopping time with respect to $\{Y_i\}$. \\

Let $P_{\omega}(m)$ be the probability that the final epoch lasts at least
$m$ spins. Let $m = kq+r, 1 \le r \le k$. Then $P_{\omega}(m) \le 
(3/4)^q$. \\

Observe that if $|S_T| \ge kn + s, \, s \ge 0$, then the last payoff $Y_T$
must satisfy $|Y_T| \ge s+k$, which is possible only if the epoch lasts at
least $s+1$ spins, since the number of tokens in the pot can go up only by
one unit at a time. Thus, $P(|S_T| \ge kn + s) \le P(|Y_T| \ge s+k) \le
P_{\omega}(s+1)$. \\

Let $s+1=kq+r, 1 \le r \le k$, so that $P_{\omega}(s+1) \le (3/4)^{q}$.
Thus, $E(|S_T|) = \sum_{i=1}^{\infty} P(|S_T| \ge i) \le kn +
\sum_{s=1}^{\infty} P(|S_T| \ge kn+s) \le kn + k \sum_{q=n}^{\infty}
(3/4)^{q}$. Similarly, $E(S_T^2) = \sum_{i=1}^{\infty}(2i-1)P(|S_T| \ge i)
\le \frac{k^2n^2}{4} + 2k \sum_{q=n}^{\infty} (q+1) (3/4)^{q}$.
Therefore, $E(|S_T|)$ is $O(n)$ and $E(S_T^2)$ is $O(n^2)$. \\

Define $n$ epochs to be an {\em age}. Note that an epoch with $k-1$ Shtels
followed by a Ganz gives $P_k$ a payoff of $2k-2$ units, and occurs with
probability $4^{-k}$. Thus the probability that all the epochs in a given
age is of the above type is $\delta \doteq 4^{-kn}$. If we ever have such
an age in a game, $P_k$ clearly wins, and we say that $P_k$ won by a {\em
landslide}. Clearly, the expected number of ages before $P_k$ wins by a
landslide is given by $\sum_{j=1}^{\infty} j(1-\delta)^{j-1}\delta =
4^{kn}$. Thus the expected number of epochs in a game of metaslowdel is at
most $n\,4^{kn}$. Therefore, $E(T)$ is finite. Similarly, it can be shown
that $E(T^2)$ is also finite. \\

Let $t=ku+v, 1 \le v \le k$. Then, $P(|Y_1| \ge t) \le
P_\omega(ku) \le (3/4)^{u-1}$. Now $|\mu| \le E(|Y_1|) =
\sum_{i=1}^{\infty} P(|Y_1| \ge i) \le k + k \sum_{u=1}^{\infty}
(3/4)^{u-1} = 5k$. \\

Suppose $|\mu| \ge \frac{1}{10}$. By Wald's equation, we have $|\mu| E(T) 
= |E(S_T)| \le E(|S_T|)$. Since $|\mu|$ is finite, and bounded below by a
positive constant, it follows that $E(T)$ is $O(n)$. \\

Now we consider the case when $|\mu| < \frac{1}{10}$. \\

Let $X$ be the collection of all sequences which form an epoch. Let $X_S$
(respectively $X_N, X_H, X_G$) consist of all sequences in $X$ whose
penultimate term is S (respectively N, H, G). For any sequence $x$ in
$X_S$, define its neighbours in $X_N, X_H, X_G$ to be the sequences which
agree with $x$ everywhere except in the penultimate position. Note that a
sequence belongs to $X_S, X_N, X_H$ or $X_G$ with probability
$\frac{1}{4}$. Clearly, a sequence in $X_S$ and its neighbour in $X_N$
cannot both have zero payoff, therefore one of them must contribute at
least one unit towards $E(Y_1^2)$. Thus, $E(Y_1^2) \ge \frac{1}{4}$.  
Therefore, $\sigma^2 = E(Y_1^2) - \mu^2 > \frac{6}{25}$. \\

But we also have, $\sigma^2 \le E(Y_1^2) \le k^2 + \sum_{t=k+1}^{\infty}
(2t-1) P(|Y_1| \ge t) \le 41k^2$. \\

By Wald's equation, we have $E(S_T^2) = \sigma^2 E(T) + \mu^2 E(T^2) \ge
\sigma^2 E(T)$. Since $\sigma^2$ is finite, and bounded below by a
positive constant, it follows that $E(T)$ is $O(n^2)$. \\

Let $p_s$ denote the probability that $P_k$ goes home after exactly $ks$
spins, and let $E_s$ denote the expected number of epochs in a game
lasting exactly $ks$ spins before $P_k$ goes home. We have, $E(T) = 
\sum_{s=1}^{\infty} p_s E_s$. \\

Given $k$, choose $\alpha$ such that $(\frac{\alpha}{64^k})^{\alpha} (1 - 
\alpha)^{1 - \alpha} > 0.76$ Let $T_{s} \doteq \lfloor \alpha s \rfloor 
< \frac{s}{2}$. Note that for sufficiently large $s$, $T_s < 
(\frac{0.76}{0.75})^{s}$. \\ 

Let $M$ denote the number of games which last $ks$ spins, with less than
$T_{s}$ epochs. We have, $M \le 4^{s(k-1)} \sum_{r=0}^{T_{s}-1} {s \choose
r} 3^{s-r}$. Since ${n \choose r} \alpha^{\alpha} (1 - \alpha)^{1 - 
\alpha} < 1$ for $r < \alpha n$, we have

\addtolength{\parskip}{-2mm}

$$M \le 4^{s(k-1)} T_{s} 3^{s} {s \choose T_{s}} <
4^{ks} {\left ( \frac{0.76}{\alpha^{\alpha} (1 - \alpha)^{(1 - \alpha)}} 
\right )}^s < 4^{ks(1 - 3 \alpha)}$$ 

\addtolength{\parskip}{3mm}

We now construct more than $M$ metaslowdel games which last $ks$ spins and
have at least $T_{s}$ epochs. Our games evolve in phases. \\

\addtolength{\parskip}{-1mm}

The first phase is restorative, (metaslowdel can start from any
configuration) and ends when there are $k$ tokens in the pot, the
difference between the number of tokens in the possession of any pair of
players is at most one, and it is the first player's turn to spin. This is
accomplished as follows: \\

We begin with a sequence of Halbs, until there are only two tokens left in
the pot. If there was only one token to begin with, we have a Shtel
instead. We then have a (possibly empty) sequence of Nishts, till it is
the first player's turn to spin. In every subsequent round, a player with
the highest number of tokens gets a Shtel, a player with the lowest number
of tokens gets a Halb, and everyone else gets Nishts. If at the end of any
round the difference between the highest and lowest is at most one, we
have $k-2$ rounds comprising a Shtel for one of the (current) leaders and
Nishts for everyone else, thus increasing the pot size from $2$ to $k$. It
is easy to see that the number of spins in the restorative phase is
$O(n)$. At the end of this phase, each player has $m + \epsilon$ tokens,
with $m \ge n-1$ and $\epsilon \in \{0,1\}$. \\

In the second phase, we have $T_{s}$ rounds in which each player gets
Ganz. This ensures that all the games we construct have at least $T_s$
epochs. \\

The third phase is divided into {\em gamelets}. A gamelet of length $\ell$
is a segment of $\ell$ spins. Note that all gamelets of length up to $n$
which start from the initial configuration of dreidel are legal, since the
payoff never drops below $-n$ or goes above $n$ for such gamelets. \\

Let $p = \lfloor \frac{n-k-1}{k^2} \rfloor$, so that $n-k-k^2 \le pk^2 <
n-k$ and let $X$ be the collection of all gamelets of length $pk+1$ that
end with a Ganz. Let $g \in X$ and define $\rho(g)=(u_1,\ldots,u_{k-1})$
if and only if $g$ gives payoffs $u_1,\ldots,u_{k-1},
-(u_1+\cdots+u_{k-1})$ for players $P_1, \ldots, P_{k-1}, P_{k}$
respectively. Let $x_{u_1, \ldots,u_{k-1}}$ denote the number of gamelets
$g$ in $X$ with $\rho(g)=(u_1,\ldots,u_{k-1})$. Note that 
$$\sum_{-pk-1 \le u_1,\ldots,u_{k-1} \le (k-1)(2p+1)} x_{u_1, 
\ldots,u_{k-1}} = 4^{pk}$$
Observe that if $g_1,g_2,\ldots,g_k$ are gamelets in $X$ with $\rho(g_1) =
\cdots \rho(g_k)$, then the concatenated gamelet $g_1g_2 \ldots g_k$ gives
zero payoff for every player. By Minkowski's inequality, there are at
least 
$$\sum_{-pk-1 \le u_1,\ldots,u_{k-1} \le (k-1)(2p+1)} x_{u_1,
\ldots,u_{k-1}}^k \ge \frac{(4^{pk})^k}{(3pk)^{k-1}} >
\frac{({\frac{k}{3})}^{k-1}}{4^{k^2+k}} \frac{4^n}{n^{k-1}}$$ 
gamelets of length $pk^2+k < n$ which give zero payoff for every player. 
For sufficiently large $n$, this number exceeds $4^{n(1 - \alpha)}$. \\

The third phase proceeds in a series of such concatenated gamelets
with payoff zero until the next gamelet would increase the number of
spins beyond $k(s-m-2)$.  When this happens, we have (at most $pk$)
rounds of Nishts till the number of spins is exactly $k(s-m-2)$. \\

The fourth and final phase has $m+2$ rounds. In the first $m$ rounds,
everyone gets Shtel. In the next round, everyone who has a token gets
Shtel, and everyone else gets Nisht. In the last round, everyone gets
Nisht, except $P_k$ who gets Ganz, wins, and goes home. \\

Observe that we have constructed at least $4^{(ks(1 - 2 \alpha)-kT_{s})} >
4^{ks(1- 3 \alpha)}$ different games lasting exactly $ks$ spins, and with
at least $T_s$ epochs for all $s \ge s_0$, where $s_0$ is sufficiently
large. By our choice of $\alpha$, this exceeds the number of possible
games with less than $T_{s}$ epochs. Thus, for all $s \ge s_0$, we have
$E_s \ge \frac{\alpha s}{2}$. \\

Let $U$ denote the number of spins in a game of metaslowdel before $P_k$ 
goes home. We have,
$$E(U)=\sum_{s=1}^{\infty} p_s ks < ks_0 + \sum_{s=s_0}^{\infty} p_s ks <
ks_0 + \frac{2}{\alpha}\sum_{s=s_0}^{\infty} p_s E_s \le k s_0+
\frac{2}{\alpha} E(T)$$

Thus the expected number of spins before $P_k$ goes home is $O(n^2)$,
irrespective of the starting configuration.  If $P_k$ wins, the game is
over. Otherwise, we have a game of metaslowdel between at most $k-1$
players. Repeating the above argument, it is easy to see that the expected
number of epochs in a game of metaslowdel is $O(n^2)$. It follows that the
expected number of spins in a game of dreidel between $k$ players is
$O(n^2)$. \\

We now present an alternative proof for the case $k=2$ using Markov chains. \\

Let $\Lambda = 2n+3$. We consider a Markov chain $M$ on an infinite state
space, $U$, where each state is indexed by a triple $(x,y,z)$ with $x$
denoting the number of tokens in the pot, $y$ denoting the number of
tokens in the possession of $P_1$ modulo $\Lambda$, and $z=i$ if and only
if $P_i$ plays next. Let $z^{*} \doteq 3-z$.  \\

The states reachable from $s_1 = (x_1,y_1,z_1)$ via a single transition
are $ \lambda_G(s_1) \doteq (2,y_1+x_1-1,z^{*}_1), \lambda_H(s_1) \doteq (\lceil
\frac{x_1}{2} \rceil, y_1+ \lfloor \frac{x_1}{2} \rfloor,z^{*}_2), \lambda_N
\doteq (x_1,y_1,z^{*}_1)$ and $\lambda_S \doteq (x_1+1,y_1-1,z^{*}_1)$. \\

The initial state is $s_0 \doteq (2,N-1,1)$. The end-states are precisely
the non-dreidel states. Finally, we partition $U$ into disjoint subsets
$A_k \doteq \{(x,y,z) \in U : x=k \}$.  \\

Observe that the first co-ordinates of the state space form a Markov chain
$M_1$. To see that the chain is irreducible, consider 
$$(x) \stackrel{\lambda_H \lambda_H \cdots \lambda_H} {\longrightarrow} (1)
\stackrel{\lambda_S \lambda_S \cdots \lambda_S} {\longrightarrow} (y)$$ 
Since the set of timesteps on which any state can be reached is cofinite 
(consider $\lambda_N \lambda_N \cdots \lambda_N$), the chain is aperiodic. Finally, 
observe that the mean return time to the state (2) is at most $E_g$, the 
expected time for a Ganz. Since $E_g = \sum_{k=1}^{\infty} k 
(\frac{3}{4})^{k-1} \frac{1}{4} = 4$, it follows that the chain is 
positive recurrent. Therefore, $M_1$ is ergodic. \\

Let $\pi_{ij}$ denote the transition probability from state $(i)$ to state
$(j)$. Let $\pi_j$ denote the stationary probability of $M_1$ being in
state $(j)$. Then, $\pi_j = \sum_{i=1}^{\infty} \pi_{ij} \pi_i$. From
these equations, it can be easily shown that $\pi_2 \ge \frac{6}{13}$, an
improvement over the easy estimate $\pi_2 \ge \frac{1}{4}$. \\

To show that $M$ is irreducible, consider an arbitrary pair of states $s_1
= (x_1,y_1,z_1), s_2 = (x_2,y_2,z_2)$. Since $\lambda_N (x_1,y_1,z_1) = 
(x_1,y_1,z^{*}_1)$, we can assume that $z_1 = z_2 = 1$. Note that 
$$(x_1,y_1,1) \stackrel{\lambda_S \lambda_N \cdots \lambda_S \lambda_N} 
{\longrightarrow} (x_1,y_2,1) \stackrel{\lambda_N \lambda_H \cdots
\lambda_N \lambda_H} {\longrightarrow} (1,y_2,1) \stackrel{\lambda_N \lambda_S \cdots
\lambda_N \lambda_S} {\longrightarrow} (x_2,y_2,1)$$

\addtolength{\parskip}{2mm}

Let $p_{ij}^k$ denote the transition probability from state $i$ to state
$j$ in exactly $k$ steps. Recall that the period of a state is the largest
integer $d$ such that $p_{ii}^n \neq 0 \Rightarrow d | n$. Since $p_{ii}^2
\ge \frac{1}{16}$ (consider $\lambda_N \lambda_N$) and $p_{ii}^{2k+1} = 0$
(consider the third co-ordinate), it follows that $d=2$. Thus all states
have period $2$. \\

\addtolength{\parskip}{-2mm}

Now consider a new Markov chain $M'$ with state space consisting of the
states $(x,y,1)$ and transition probabilities given by $q_{ij} =
p_{ij}^2$, where $p_{ij}$ are the transition probabilities of $M$. It
follows from the above arguments that $M'$ is irreducible and aperiodic.
Note that for any fixed $i, \sum_{j \in A_2} q_{ij}^n \ge \frac{1}{4}$ for
all $n$. Since $|A_2|$ is finite, there exists $j \in A_2$, such that
$\lim_{n \rightarrow \infty} q_{ij}^n > 0$. Thus there exist stationary
probabilities $\pi_j$. \\

A sequence starting at the initial state is said to be {\em fast} if it
reaches an endstate before returning to the initial state, and is said to
be {\em slow} if it reaches the initial state before returning to the end
state. Let $p_f$ (respectively $p_s$) denote the probability that a
sequence starting at the initial state is a fast (respectively slow)  
sequence. Since the chain is positive recurrent, the sequence returns to
the initial state with probability $1$. Therefore, $p_f+p_s=1$. \\

Let $\mu_0$ denote the mean return time, i.e., the expected number of
steps to return to the initial state. We have, $\mu_0= p_f \mu_f + p_s
\mu_s$, where $\mu_f$ and $\mu_s$ denote the mean return times for fast
and slow sequences. \\

We note that the definition of the second co-ordinate ensures that it is
not possible to make an illegal move from a dreidel state to another
dreidel state without passing through an end state. Therefore, a dreidel
game ends without returning to the initial state with probability $p_f$
and returns to the initial state before ending with probability $p_s$. The
former shall be called fast games and the latter, slow games. \\

Let $\mu_d$ denote the mean duration of a dreidel game, and let $\mu_{df}$
and $\mu_{ds}$ denote the mean duration of fast and slow dreidel games
respectively. Note that $\mu_{ds} = \mu_s + \mu_d$ and $\mu_{df} \le
\mu_{f}$. Since $\mu_{d}=p_f \mu_{df} + p_s \mu_{ds}=p_s \mu_{df} +p_s
(\mu_s + \mu_d)$, it follows that $$\mu_d = \frac{p_f \mu_{df} + p_s
\mu_{ds}}{1-p_s} \le \frac{p_f \mu_f + p_s \mu_s}{p_f} =
\frac{\mu_0}{p_f}$$

Since $\pi_j = \pi_k$ for all $j, k \in A_2$ (by symmetry), we have $\pi_j
= \frac{\pi_2}{\Lambda}$. Let $\mu'_0$ denote the mean return time for the
initial state in $M'$. We have, $\mu'_0 = \frac{1}{\pi_j} =
\frac{\Lambda}{\pi_2}$. It follows that $\mu_0 =\frac{2 \Lambda}{\pi_2}
\le \frac{13(2N+3)}{3}$ \\

We now derive a lower bound for $p_f$. \\

Let $P[y_1,z_1;y_2,z_2;y_2,z_3]$ denote the probability of reaching
$(2,y_2,z_2)$ before $(2,y_3,z_3)$ given that we start at $(2,y_1,z_1)$.
By an extension of notation, given a set of states $S$ in $M$,
$P[y_1,z_1;y_2,z_2;S]$ shall denote the probability of reaching
$(2,y_2,z_2)$ before any of the states in $S$ given that we start at
$(2,y_1,z_1)$. \\

Let $\overline{a}, a \oplus b$ and $a \ominus b$ denote $-a$ mod
$\Lambda$, $a+b$ mod $\Lambda$ and $a-b$ mod $\Lambda$ respectively. The
following identities are easily verified: \\

\begin{itemize}

\item {\em Duality:} $P[y_1,z_1;y_2,z_2;y_3,z_3] = P[\overline{y_1} 
\ominus 2, z^{*}_1; \overline{y_2} \ominus 2, z^{*}_2; \overline{y_3} 
\ominus 2, z^{*}_3]$

\item {\em Complementarity:} $P[y_1,z_1;y_2,z_2;y_3,z_3] = 1 - 
Pr[y_1,z_1;y_3,z_3;y_2,z_2]$

\item {\em Translation Invariance:} $P[y_1,z_1;y_2,z_2;y_3,z_3] = P[y_1 
\oplus m, z_1; y_2 \oplus m, z_2; y_3 \oplus m, z_3]$

\end{itemize} 

\vspace*{4mm} 

Let $A^{y_1}_m = P[y_1,1;y_1 \oplus m, 2; y_1 \ominus 1,1]$. We have, \\

$\frac{A^{y_1}_{m+1}}{A^{y_1}_{m}} 
\ge P[y_1 \oplus m,2;y_1 \oplus (m+1),2;y_1 \ominus 1,1]$ \\
 
$= P[\overline{(y_1 \oplus m)} \ominus 2,1;(\overline{y_1 \oplus (m+1)}) 
\ominus 2, 1; (\overline{y_1 \ominus 1}) \ominus 2,2]$ ({\em {Duality}}) \\

$= P[y_1,1; y_1 \ominus 1, 1; y_1 \oplus (m+1),2]$ ({\em {Translation 
Invariance}}) \\

$= 1 - P[y_1,1; y_1 \oplus(m+1),2; y_1 \ominus 1, 1]$ ({\em 
{Complementarity}}) \\
 
$= 1 - A^{y_1}_{m+1}$ \\

Since $A^{y_1}_1 \ge 1/4$ (consider $\lambda_G$), it follows from induction 
that $A^{y_1}_m \ge \frac{1}{m+3}$. \\

Let $B^{y_1}_m \doteq Pr[y_1,1;y_1 \ominus m, 2; y_1 \oplus 1,1]$.  As
before, it can be shown that $\frac{B^{y_1}_{m+1}}{B^{y_1}_{m}} \ge 1 -
B^{y_1}_{m+1}$. Since $B^{y_1}_1 \ge 1/64$ (consider $\lambda_S \lambda_H
\lambda_N$), it follows from induction that $B^{y_1}_m \ge \frac{1}{m+63}$. \\

Let $S_1 = \{(2,N-1,1),(2,N-2,1) \}$ and $S_2 = \{(2,N-1,1),(2,N,1) \}$.  
Note that $\omega_1 \doteq P[N-1,1;N,1;S_1] \ge \frac{1}{8}$ (consider
$\lambda_G \lambda_N$ and $\lambda_H \lambda_S$). Similarly, $\omega_2 \doteq
P[N-1,1;N-2,1;S_2] \ge \frac{1}{8}$ (consider $\lambda_N \lambda_G$ and $\lambda_S
\lambda_H$). Now, \\

$p_f \ge P[N-1,1; 2N+1,2; N-1,1]$ \\

$\ge \omega_1 P[N,1; 2N+1,2; N-1,1] + \omega_2 P[N-2,1;2N+1,2;N-1,1]$ \\
 
$\ge \frac{1}{8} A^{N}_{N+1} + \frac{1}{8} B^{N-2}_N \ge \frac{1}{8(N+4)} + 
\frac{1}{8(N+63)} \ge \frac{1}{4(N+63)}$ \\

Thus, $\mu_d \le \frac{\mu_0}{p_f} \le \frac{104N^2}{3} +o(N^2)$.  
This completes the proof. \\

\newpage

{\noindent {\bf Acknowledgements}} \\

We thank Professors J\'{o}zsef Beck, J\'{a}nos Koml\'{o}s and Doron
Zeilberger for insightful suggestions and several useful discussions.
Thanks are also due to the attendees of the Graduate Student Combinatorics
Seminar at Rutgers, most notably Stephen Hartke, Vincent Vatter and
Nicholas Weininger, for their careful scrutiny of the first version of our
proof. \\

{\noindent {\bf References}}

\begin{enumerate}

\item E. Demaine, R. Fleischer, A. S. Fraenkel, R.J. Nowakowski, {\em 
{Open Problems at the 2002 Dagstuhl Seminar on Algorithmic Combinatorial 
Game Theory}}, Theoretical Computer Science $313$, $2004$.

\item R. Feinerman, {\em {An Ancient Unfair Game}}, The American 
Mathematical Monthly $83$, $1976$.

\item F. M. Trachtenberg, {\em {The Game of Dreidel Made Fair}}, The 
College Mathematics Journal $27$, $1996$.

\item D. Zeilberger. {\em Does Dreidel Last $O(NUTS^2)$ Spins?}

\hspace*{-1.0cm} 
http://www.math.rutgers.edu/~zeilberg/mamarim/mamarimhtml/dreidel.html

\end{enumerate}

\end{document}